\documentclass[12pt]{article}

\usepackage{arxiv}

\usepackage{fullpage}
\usepackage[T1]{fontenc}
\usepackage{microtype}
\usepackage[authoryear,round]{natbib}

\usepackage[colorlinks=true, citecolor=cyan, linkcolor=cyan, urlcolor=cyan]{hyperref}
\usepackage{orcidlink}

\usepackage{graphicx}
\usepackage{amsmath,amssymb,amsfonts,mathtools}
\usepackage{multirow}
\usepackage{booktabs}
\usepackage{makecell}
\usepackage{tabularx}
\usepackage{algorithm}
\usepackage{algpseudocode}
\usepackage{xcolor}
\usepackage{tikz}
\usepackage{pgfplots}
\usepgfplotslibrary{groupplots}
\usepgfplotslibrary{statistics}
\pgfplotsset{compat=1.18}
\usepackage{cleveref}
\usepackage{eurosym}

\definecolor{mycol1}{HTML}{7CA3CB} 
\definecolor{mycol2}{HTML}{DDCC77} 
\definecolor{mycol3}{HTML}{CC6677} 
\definecolor{mycol4}{HTML}{117733} 
\definecolor{mycol5}{HTML}{882255} 
\definecolor{mycol6}{HTML}{44AA99} 

\newcommand{\reals}{\mathbb{R}}

\graphicspath{{figures/}}

\title{Emissions and electricity costs of residential PV-battery systems under self-consumption and predictive control}
\author{\mdseries
	Anna Sacilotto\thanks{
		Department of Mathematics ``F. Casorati'',
		University of Pavia, 27100 Pavia, Italy.\\
		\texttt{anna.sacilotto01@universitadipavia.it}.
	}
	\and
	Ambrogio Maria Bernardelli \orcidlink{0000-0002-2328-7062}\thanks{
		Department of Mathematics ``F. Casorati'',
		University of Pavia, 27100 Pavia, Italy.
		\texttt{ambrogiomaria.bernardelli@unipv.it}. Corresponding author.
	}
	\and
	Stefano Gualandi \orcidlink{0000-0002-2111-3528}\thanks{
		Department of Mathematics ``F. Casorati'',
		University of Pavia, 27100 Pavia, Italy.
		\texttt{stefano.gualandi@unipv.it}.
	}
}

\begin{document}
\maketitle

\begin{abstract}
Residential batteries commonly follow a self-consumption rule that ignores when grid electricity is clean or cheap.
We assess when optimization improves on it, using measured loads of 794 Swiss households, satellite-based PV at three Italian sites, archived weather forecasts and Italian grid data.
We compute the gap between the rule and perfect foresight for imported-electricity emissions, with average and marginal factors, and for retail cost.
On 20 selected cases we simulate model predictive control and a hybrid policy in which the rule executes setpoints from an hourly linear program.
In the base configuration, with average factors the rule leaves a median of 6\% of the attainable benefit in houses and under 1\% in apartments, and forecast-based policies do worse than it in 7 of 8 cases.
With marginal factors the median gap of houses is 23\%; with cost it is 36 euros for houses and 41 for apartments.
Scenario-based control recovers 76--88\% of the marginal-emission gap in three houses (70--74\% with previous-year factors) and 47--73\% of the cost gap in six cases; the hybrid policy recovers 70--85\% and 34--68\%.
Current measurements matter more than weather forecasts.
Crediting exports at the marginal factor changes the emission comparison.

\end{abstract}


\section{Introduction}\label{sec:intro}

Residential photovoltaic (PV) systems are increasingly installed together with a battery energy storage system (BESS), which shifts the PV surplus of the day to the evening and night.
A common operating strategy is the self-consumption rule: the battery charges with the PV surplus and discharges to cover the deficit, using only the measurements of the current instant~\citep{luthander2015photovoltaic,azuatalam2019energy}.
The rule needs no forecasts, never charges from the grid, and ignores when grid electricity is cheap or clean.
A controller that anticipates load, PV generation and the state of the grid could do better, and a large literature proposes optimization-based controllers for this purpose.
This paper quantifies how much a household can gain, in carbon emissions or in cost, by replacing the rule with an optimization-based policy, and how much of that gain predictive policies obtain with the forecasts available in practice.

\paragraph{Emissions and cost.}
The environmental benefit of a residential battery is not guaranteed.
In 16 utility territories of the United States, residential storage operated to minimize the electricity bill mostly increased emissions, whatever its operating mode~\citep{babacan2018unintended}.
Across 22 regions, cost-only scheduling of a household battery cut energy costs by $34\%$ on average but raised emissions by $7.1\%$ on average, with increases in 20 regions~\citep{olivieri2020optimization}.
For behind-the-meter batteries of commercial and industrial buildings, net emissions are driven primarily by conversion losses rather than by the timing of charging and discharging~\citep{fisher2017emissions}.
Operating a battery for emissions therefore requires a carbon signal of the grid and the right one: marginal emission factors, rather than average ones, measure the emissions caused or avoided by a change in demand~\citep{hawkes2010marginal,silerevans2012marginal}, and average factors can misestimate them by about $30$--$35\%$.
Recent studies address the sizing of stand-alone residential batteries for cost and emission savings~\citep{aikman2026sizing} and the low-carbon operation of plug-in batteries in residential buildings~\citep{amaripadath2026achieving}; here the PV-battery system is given, and the question is how to operate it.
We consider both emission factors, together with the retail cost of electricity.

\paragraph{Battery control as a sequential decision problem.}
Operating a battery is a sequential decision problem: at each period the controller observes the state of charge (SOC) and the realized load and generation, decides how much to charge or discharge, and then observes new information.
The object of design is a \emph{policy}, a rule that maps the available information to a decision~\citep{powell2016tutorial,powell2019unified}.
The self-consumption rule is a policy function approximation in the taxonomy of~\citet{powell2016tutorial}; model predictive control (MPC), deterministic~\citep{parisio2014model,rawlings2017model} or scenario-based~\citep{mesbah2016stochastic,vandermeer2021alternative}, is a direct lookahead approximation that solves an optimization model over a finite horizon at every period and implements the first decision.
Model predictive control is the most common way to handle forecast errors in home energy management, and it is usually evaluated against a perfect-knowledge counterpart~\citep{beaudin2015home}.
Reviews of the field also observe that results are hard to compare, for lack of a common baseline~\citep{beaudin2015home}, and that more sophisticated strategies do not necessarily perform better than the self-consumption rule when forecasts are imperfect~\citep{azuatalam2019energy}.

\paragraph{Planners and real-time rules.}
A practical way to combine anticipation with the robustness of the rule is a two-layer controller, in which an optimization layer with forecasts guides a fast layer that acts on the measurements.
The upper layer has been a day-ahead dynamic program~\citep{riffonneau2011optimal}, a linear program solved every hour~\citep{hanna2014energy}, a convex program solved at the start of each tariff period~\citep{wang2014adaptive}, a stochastic program solved every 15 minutes~\citep{hafiz2020realtime}, a two-stage program~\citep{langenmayr2020unit}, or a mixed-integer MPC~\citep{elkazaz2020hierarchical}.
What is passed to the lower layer varies: a power schedule~\citep{riffonneau2011optimal,langenmayr2020unit} or set point~\citep{hafiz2020realtime}, a threshold on the power drawn from the grid~\citep{hanna2014energy,wang2014adaptive}, a static energy reserve~\citep{langenmayr2020unit}, or the planned trajectory of the stored energy~\citep{elkazaz2020hierarchical}, which the rule-based layer uses as a discharge floor during the day and as a grid-charging target at night.
Closer to the self-consumption rule, \citet{kazhamiaka2016practical} modify it with a single daily SOC threshold, set from a forecast or from past data, that acts as charging target and discharge floor in the low-price hours, and report that such simple rules are near-optimal; forecast-based strategies for home storage limit the SOC to the energy expected to be needed during the night, mainly to extend battery life~\citep{angenendt2018comparison}.

\paragraph{This paper.}
Earlier work thus shows that prediction and hierarchical control can improve on the rule, and that simple rules are often hard to beat.
What is missing is a joint picture, under a common baseline and for different objectives, of how much improvement is attainable at all and how much of it a deployable controller obtains.
We use measured household loads, PV generation computed from satellite irradiance at Italian sites together with archived weather forecasts, and Italian grid data of 2023--2024, and we address three questions.
First, how much can any policy gain over the self-consumption rule, for which households, system sizes and objectives?
Second, how much of that gain do forecast-based MPC policies recover, and how does it depend on the weather forecasts and on measuring the current period?
Third, can the gain be obtained by letting the rule execute a plan, without real-time optimization?

\paragraph{Contributions.}
\begin{itemize}
    \item A screening assessment of the headroom of optimization, i.e.\ the gap between the self-consumption rule and a perfect-foresight benchmark that bounds the benefit of any policy, for $794$ measured household loads, three Italian sites, three PV and three battery sizes, and three objectives (average and marginal emission factors, and retail cost), with sensitivity to the accounting of exports and to the battery assumptions.
    \item A closed-loop comparison, on $20$ selected cases over $318$ days of 2024 with the state of charge carried over from one period to the next, of deterministic and scenario-based MPC under explicit assumptions on the availability of weather forecasts and grid signals, with variants that separate the role of weather forecasts, of measuring the current period, and of perfect information over the lookahead window, and a sensitivity in which the marginal emission factors are taken from the previous year.
    \item A hybrid policy that keeps the self-consumption rule as the real-time layer and translates the plan of a lookahead linear program into four SOC setpoints (floor, ceiling, target and dump); with no setpoints it is exactly the rule. The interface combines the floor, target and ceiling mechanisms used separately in earlier work~\citep{kazhamiaka2016practical,elkazaz2020hierarchical,angenendt2018comparison} and adds an export setpoint.
\end{itemize}
The input data are public, and the code and processed data are available at~\citet{bernardelligithub}.

\paragraph{Outline.}
\Cref{sec:problem} states the problem and the sequential setting, \cref{sec:policies} the policies, and \cref{sec:scenarios} the scenario model, whose details are in Appendix~\ref{app:scenarios}.
\Cref{sec:data} describes the data, \cref{sec:headroom,sec:closed-loop} the experiments, and \cref{sec:discussion,sec:conclusions} discuss the results and conclude.

\section{Problem description}\label{sec:problem}

\subsection{System and objectives}\label{sec:system}

We consider a household with PV generation and a battery, connected to the grid.
Time is divided into periods $t = 1, 2, \ldots$ of length $\Delta t = 15$~minutes, grouped into days of $T = 96$ periods.
In period $t$, let $l_t$ be the load, $s_t$ the PV generation, $b^+_t$ and $b^-_t$ the discharging and charging power of the battery, and $y^+_t$ and $y^-_t$ the power imported from and exported to the grid, all in kW.
The power balance is
\begin{equation}
    y^+_t - y^-_t = l_t - s_t - b^+_t + b^-_t,
    \label{eq:balance}
\end{equation}
and the energy $B_t$ stored at the beginning of period $t$ evolves as
\begin{equation}
    B_{t+1} = B_t + \Delta t \left( \eta b^-_t - b^+_t / \mu \right),
    \label{eq:dynamics}
\end{equation}
where $\eta$ and $\mu$ are the charging and discharging efficiencies.
The battery has usable capacity $\bar{B}$ and power $\bar{P}_b$, so that $0 \le B_t \le \bar{B}$, $0 \le b^-_t \le \bar{P}_b$ and $0 \le b^+_t \le \mu \bar{P}_b$: $b^\pm_t$ are measured on the household side, and the power drawn from the cells, $b^+_t / \mu$, is limited to $\bar{P}_b$.
We refer to the stored energy $B_t$ (kWh) as the state of charge (SOC).
To limit degradation, the energy exchanged by the battery in each calendar day $d$ is limited to $c$ equivalent full cycles,
\begin{equation}
    \Delta t \sum_{t \in d} \left( b^+_t / \mu + \eta b^-_t \right) \le 2 c \bar{B}.
    \label{eq:cycles}
\end{equation}
We do not model limits on the power exchanged with the grid.

The objective is the sum over the periods of
\begin{equation}
    \Delta t \left( p^+_t y^+_t - p^-_t y^-_t + w\, b^+_t \right),
    \label{eq:objective}
\end{equation}
where $p^+_t$ and $p^-_t$ are the value of a kWh imported and exported in period $t$, and $w$ is a cost per kWh delivered by the battery.
Three instances are considered.
With \emph{average emission factors}, $p^+_t = \alpha_t$ is the average carbon intensity of the national generation mix (gCO$_2$/kWh), exports are not credited ($p^-_t = 0$), and $w = \beta$ is the life-cycle emission factor of the energy delivered by the battery.
With \emph{marginal emission factors}, $p^+_t = \alpha^\mathrm{m}_t$ is the marginal emission factor of the grid, with $p^-_t = 0$ and $w = \beta$.
Both emission objectives count grid imports and the allocated battery term; neither credits the displacement of grid generation by exports. They therefore evaluate an import-based emission metric, rather than the net change in electricity-system emissions.
With \emph{cost}, $p^+_t$ is the retail price of electricity, $p^-_t$ the price paid for the exported energy, and $w$ a wear cost of the battery (\cref{sec:data}).
In all cases $0 \le p^-_t \le p^+_t$ and $w \ge 0$; together with the conversion losses, this implies that simultaneous import and export, or simultaneous charging and discharging, cannot reduce the objective, so that no binary variables are needed.

When a policy is evaluated over a horizon of $N$ periods, a final SOC below the initial one is charged the cost of restoring the deficit from the grid, $(B_1 - B_{N+1})^+ \, \bar{p}^+ / \eta$, where $\bar{p}^+$ is the average of $p^+_t$ over the horizon; no credit is given for energy left in the battery.
This is an evaluation convention common to all policies, not a physical restoration schedule.

\subsection{Sequential decisions and information}\label{sec:sequential}

The load, the PV generation and, for some objectives, the coefficient $p^+_t$ are uncertain and revealed progressively.
Let $\xi_t = (l_t, s_t, p^+_t)$.
A policy maps the information available at decision time to the battery decisions $(b^+_t, b^-_t)$; the grid exchanges then follow from the balance~\eqref{eq:balance}, so that any mismatch between the battery decision and the realized net load $l_t - s_t$ is absorbed by the grid.
We distinguish two information structures.
With \emph{real-time} information, the realization $\xi_t$ of the current period is known when $(b^+_t, b^-_t)$ is decided.
This idealizes, at the resolution of the data, a controller that measures the instantaneous imbalance between generation and load and adjusts the battery continuously; what happens within a period is not modelled.
With \emph{delayed} information, the decision is based on $\xi_1, \ldots, \xi_{t-1}$ only, as for a controller that sends a set point before the period starts.
Both structures use the current SOC $B_t$.
The quality of a policy is the realized value of the objective over the simulation horizon.

\subsection{Self-consumption rule}\label{sec:rule}

\begin{algorithm}[t]
\caption{Self-consumption rule with setpoints. With $f_t = g_t = 0$ and $u_t = v_t = \infty$ it is the plain self-consumption rule (SC).}
\label{alg:rule}
\begin{algorithmic}[1]
\Require SOC $B_t$, measured net load $n_t = l_t - s_t$, residual cycle budget $R_t$, setpoints $f_t$ (floor), $u_t$ (ceiling), $g_t$ (target), $v_t$ (dump)
\State $\bar{c} \gets \min\{\bar{P}_b,\ (\bar{B} - B_t)/(\eta \Delta t),\ R_t/(\eta \Delta t)\}$ \Comment{charging limit}
\State $b^-_t \gets \min\{(-n_t)^+,\ \bar{c},\ (u_t - B_t)^+/(\eta \Delta t)\}$ \Comment{store the PV surplus up to the ceiling}
\If{$B_t < g_t$} $b^-_t \gets \max\{b^-_t,\ \min\{(g_t - B_t)/(\eta \Delta t),\ \bar{c}\}\}$ \Comment{charge towards the target}
\EndIf
\If{$b^-_t > 0$} $b^+_t \gets 0$
\Else
    \State $\bar{d} \gets \min\{\mu \bar{P}_b,\ \mu R_t / \Delta t\}$ \Comment{discharging limit}
    \State $b^+_t \gets \min\{n_t^+,\ \mu (B_t - f_t)^+ / \Delta t,\ \bar{d}\}$ \Comment{cover the deficit down to the floor}
    \If{$B_t > v_t$} $b^+_t \gets \max\{b^+_t,\ \min\{\mu (B_t - v_t)/\Delta t,\ \bar{d}\}\}$ \Comment{export down to the dump level}
    \EndIf
\EndIf
\end{algorithmic}
\end{algorithm}

The self-consumption rule (SC) is a common residential operating policy~\citep{luthander2015photovoltaic,azuatalam2019energy}.
At each period it charges the battery with the PV surplus, within the limits of capacity, power and residual cycle budget, and exports the rest; when PV does not cover the load, it discharges the battery to cover the deficit and imports the rest.
It uses real-time information, needs no forecast, never charges from the grid, and ignores $p^\pm_t$.
\Cref{alg:rule} states it in a general form, with four setpoints that will be used by the hybrid policies of \cref{sec:hybrid}: the battery does not discharge below the floor $f_t$ to cover the load, does not charge from the PV surplus above the ceiling $u_t$, charges (from the grid if needed) up to the target $g_t$, and discharges towards the dump level $v_t$ even beyond the load, exporting the excess.
With $f_t = g_t = 0$ and $u_t = v_t = \infty$, \cref{alg:rule} is exactly SC.

\section{Optimization-based policies}\label{sec:policies}

All optimization-based policies solve linear programs over a lookahead window of $H = 96$ periods (24 hours) that starts at the current period $\tau$ and therefore covers the rest of the current day and the beginning of the next one.
The window is described by $M$ scenarios $\xi_{t,k}$, $t = \tau, \ldots, \tau + H - 1$, $k = 1, \ldots, M$, generated by the model of \cref{sec:scenarios}.

\subsection{Lookahead models}\label{sec:lookahead}

The \emph{deterministic} lookahead model uses a single trajectory $(\hat{l}_t, \hat{s}_t, \hat{p}^\pm_t)$, the average of the scenarios, and is the linear program
\begin{equation}\label{eq:lookahead}
    \begin{array}{ll}
        \text{minimize} & \Delta t \sum_{t=\tau}^{\tau+H-1} \left( \hat{p}^+_t y^+_t - \hat{p}^-_t y^-_t + w\, b^+_t \right) - \epsilon \sum_{t=\tau}^{\tau+H-1} B_{t+1}\\
        \text{subject to} & y^+_t - y^-_t = \hat{l}_t - \hat{s}_t - b^+_t + b^-_t, \quad B_{t+1} = B_t + \Delta t (\eta b^-_t - b^+_t/\mu),\\
        & 0 \le B_{t+1} \le \bar{B}, \quad 0 \le b^-_t \le \bar{P}_b, \quad 0 \le b^+_t \le \mu \bar{P}_b, \quad y^\pm_t \ge 0, \quad t = \tau, \ldots, \tau + H - 1,\\
        & \text{cycle limit~\eqref{eq:cycles}, separately for the current and for the next day},
    \end{array}
\end{equation}
where $b^\pm_t$, $y^\pm_t$ and $B_{t+1}$ are the optimization variables and the measured SOC $B_\tau$ is problem data.
The cycle limit of the current day applies to the periods of the window before midnight, with the budget reduced by the energy already exchanged since midnight; the periods after midnight have the full budget of the next day, and the budget not used today is not carried over.
No terminal condition or terminal value is imposed on the lookahead SOC; the 24-hour window is a modelling choice and does not eliminate end-of-window effects.
Because exported energy has no value under the emission objectives, the model is often indifferent between storing the PV surplus now or later; the last term, with $\epsilon$ fixed for each objective to $10^{-5}$ times the median of $p^+_t$ over the data set, favours storing early among solutions of nearly equal value, and is not part of the realized objective.

The \emph{two-stage} stochastic model uses all $M$ scenarios: the decision $(b^+_\tau, b^-_\tau)$ of the first period is common to all scenarios, while all subsequent decisions depend on the scenario, and the objective is the average over the scenarios of the objective of~\eqref{eq:lookahead}.
As usual for two-stage models, the recourse decisions assume that the scenario is revealed after the first period.
It is a linear program with $O(MH)$ variables.

\subsection{MPC policies}\label{sec:mpc}

An MPC policy solves a lookahead model at every period $\tau$ with the measured SOC, implements the first decision, projected if necessary onto the decisions that are feasible given the current SOC and cycle budget, and repeats at the next period.
We consider:
\begin{itemize}
    \item MPC-D, with the deterministic model on the forecast average, and MPC-S, with the two-stage model on $M = 100$ scenarios, both with real-time information: the first period of every scenario is the measured one;
    \item MPC-S0, as MPC-S but with scenarios that are not conditioned on the weather forecast, which shows the effect of the forecast-conditioned scenario model (the model is refitted without the forecast, so this is not a comparison under an otherwise fixed distribution);
    \item MPC-SL, as MPC-S but with delayed information, which shows the value of measuring the current period;
    \item MPC-P, the deterministic model with the actual realizations over the window, as a perfect-information reference for the chosen 24-hour lookahead formulation. Its realized objective is not a bound on all forecast-based policies.
\end{itemize}

\subsection{Hybrid policies}\label{sec:hybrid}

A hybrid policy keeps the self-consumption rule as the real-time layer and uses optimization only to set its four setpoints.
A plan is computed by solving the deterministic lookahead model~\eqref{eq:lookahead} with delayed information; let $\hat{b}^\pm_t$ and $\hat{B}_{t+1}$ be the planned flows and SOC, and $\hat{n}_t = \hat{l}_t - \hat{s}_t$ the forecast net load.
The plan is translated into setpoints for each period of its validity:
\begin{itemize}
    \item if the plan charges more than the forecast PV surplus ($\hat{b}^-_t > (-\hat{n}_t)^+$), the target is $g_t = \hat{B}_{t+1}$: the battery charges from the grid;
    \item if it charges less than the surplus ($\hat{b}^-_t < (-\hat{n}_t)^+$), the ceiling is $u_t = \hat{B}_{t+1}$: the remaining PV is exported;
    \item if it discharges less than the deficit ($\hat{b}^+_t < \hat{n}_t^+$), the floor is $f_t = \hat{B}_{t+1}$: energy is held for later;
    \item if it discharges more than the deficit ($\hat{b}^+_t > \hat{n}_t^+$), the dump level is $v_t = \hat{B}_{t+1}$: the battery exports;
\end{itemize}
and the setpoints take their neutral values ($0$ or $\infty$) otherwise.
The rule then acts on the measured net load every 15 minutes.
The setpoints are defined so that, if the realized net load coincides with the forecast, the rule reproduces the planned flows; if it does not, the rule still behaves as SC within the bounds set by the plan, so that, for example, an unexpected deficit is covered from the battery unless the plan reserved that energy for later.
We consider HYB-1h, in which the plan is recomputed every hour and executed for the next four periods, and HYB-24h, in which it is computed once a day at midnight; HYB-P recomputes the plan every hour with the actual realizations, as a perfect-information reference for this hybrid implementation.
A hybrid policy solves a single small linear program per hour (or per day) instead of a stochastic program every 15 minutes, and the battery controller executes a rule on its own measurements, as it already does with SC.

\subsection{Benchmarks}\label{sec:benchmarks}

The perfect-foresight benchmark (PF) solves a single linear program over the whole simulation horizon with the actual realizations, the same initial SOC, the daily cycle limits, and the terminal charge of \cref{sec:system}; its value is a lower bound on the realized objective of every policy.
The household without battery (NB) gives the reference for the benefit of the battery.
\Cref{tab:policies} summarizes the policies.

\begin{table}[t]
\caption{Policies. Information: RT = the current period is measured before deciding; delayed = only the previous periods; full = perfect information.}
\label{tab:policies}
\centering
{\small
\setlength{\tabcolsep}{3pt}
\begin{tabularx}{\linewidth}{@{}l>{\raggedright\arraybackslash}Xll@{}}
\toprule
Policy & Description & Information & Optimization \\
\midrule
SC & self-consumption rule (\cref{alg:rule}) & RT & none \\
MPC-D & deterministic lookahead on the forecast average & RT & every 15 min \\
MPC-S & two-stage stochastic lookahead, 100 scenarios & RT & every 15 min \\
MPC-S0 & as MPC-S, scenarios not conditioned on the weather forecast & RT & every 15 min \\
MPC-SL & as MPC-S, decision before measuring the period & delayed & every 15 min \\
MPC-P & deterministic lookahead with the actual realizations & full (24 h) & every 15 min \\
HYB-1h & SC with setpoints from a deterministic plan & RT (rule), delayed (plan) & every hour \\
HYB-24h & as HYB-1h, plan computed at midnight & RT (rule), delayed (plan) & once a day \\
HYB-P & as HYB-1h, plan with the actual realizations & full (24 h) & every hour \\
PF & perfect foresight over the whole horizon & full & once \\
\bottomrule
\end{tabularx}
}
\end{table}

\section{Scenario model}\label{sec:scenarios}

The lookahead models need scenarios of load, PV generation and, when it is not known in advance, the grid signal (carbon intensity or day-ahead price), conditioned on what has been observed so far and on the weather forecasts available at decision time.
We extend the scenario generation of~\citet{sacilotto2024riduzione}, which samples daily profiles from the leading principal components of the preceding days, to a joint model of all the signals that is conditioned on the observations and forecasts.
For each of the $28$ days preceding the current one, the daily profiles of load, PV generation, logarithm of the grid signal and PV forecast are stacked into a vector of $4T$ entries ($3T$ for MPC-S0, without the forecast); the first $K = 5$ principal components are retained, and the model is interpreted as a probabilistic principal component analysis~\citep{tipping1999probabilistic}.
Given the entries observed at decision time, the latent coefficients have a Gaussian posterior, from which the scenarios of the unobserved entries are sampled; Appendix~\ref{app:scenarios} gives the details.

The information used at period $\tau$ is the following.
Load, PV generation and grid signal of the current day are observed up to period $\tau$ (included with real-time information, excluded with delayed information).
The weather forecast is the value predicted $24$ hours before each valid time (\cref{sec:data}); allowing $6$ hours for the publication of the model run and for the interpolation of hourly values, a conservative assumption rather than a reconstruction of the publication time of each run, at period $\tau$ we use only the forecasts of the valid times up to $18$ hours ahead, which cover the rest of the current day and, in the afternoon and evening, the first hours of the next one.
The last hours of the window therefore rely on the statistical model alone.
Day-ahead prices are known for the current day and, from 14:00, for the next one.
The marginal emission factors are a table by quarter and hour (\cref{sec:data}) and are supplied as known, whereas the average factor is forecast by the scenario model: the two emission objectives therefore differ also in the information available on the signal.
MPC-S0 fits the same model without the forecast block.

\section{Data and case study}\label{sec:data}

All series cover 2023 and 2024 and are aligned on a 15-minute grid in UTC+1 throughout the year, so that every day has $96$ periods: the loads are metered every 15 minutes, the grid signals are hourly and held constant within the hour, and PV generation is computed from hourly irradiance and interpolated.
Policies are simulated in closed loop from 17 February to 30 December 2024 ($318$ days), starting from a half-full battery; the scenario model uses the $28$ days preceding each day, and complete archived forecasts are available from 20 January 2024.

\paragraph{Household loads.}
Household loads are the 15-minute smart-meter measurements of 2{,}447 residential installations of a Swiss utility in 2023 and 2024, published by~\citet{kaiser2026swissloads}.
We use them as demand trajectories for a comparative assessment under Italian weather and grid conditions: the data set provides individual measurements at 15 minutes over the same two years as the grid data, household-level variation, and equipment metadata, and it is public.
All policies face the same trajectories within each case, so the comparisons are conditional on these demand profiles and are not estimates for the Italian residential population.
To restrict the assessment to household demand without recorded electrified heating, vehicle charging or PV, we keep single-family houses and apartments without PV and without any recorded heat pump, electric water heater, electric heating, electric vehicle charger or change of equipment; we discard installations with less than $98\%$ of the values, those whose average winter consumption exceeds $1.6$ times the summer one, a screen against unrecorded electric heating, and three apartments that consumed less than $100$~kWh in 2024, which suggests that they were unoccupied.
Missing values ($0.2\%$ on average, at most $1.9\%$ per household) are interpolated over gaps of up to one hour and otherwise replaced by the value of the same period one or two weeks earlier or later.
The result is $49$ houses (median consumption $5{,}232$~kWh in 2024) and $745$ apartments (median $1{,}978$~kWh, close to the $2{,}000$~kWh per year of the average Italian resident household~\citep{arera2024schedaQ3}); this is a check of scale only, not of the daily timing, variability or cooling demand of the loads.
Using the same calendar days as the grid and weather data preserves the calendar alignment; it does not reproduce the response of Italian household demand to local weather. Climate differences are discussed in \cref{sec:discussion}.

\paragraph{PV generation and forecasts.}
We place the household in Milan, Rome or Palermo, with a south-facing array tilted by $30^\circ$.
The actual PV generation is computed from the hourly global tilted irradiance of the SARAH-3 satellite data set~\citep{pfeifroth2023sarah3} and the ERA5 air temperature~\citep{hersbach2020era5}, both retrieved from the Open-Meteo archive~\citep{zippenfenig2023openmeteo}, with the model
\begin{equation*}
    s_t = P_\mathrm{PV} \frac{G_t}{1000} \left(1 + \gamma (T^\mathrm{c}_t - 25)\right) (1 - L), \qquad T^\mathrm{c}_t = T^\mathrm{a}_t + G_t \frac{T_\mathrm{NOCT} - 20}{800},
\end{equation*}
where $G_t$ is the irradiance (W/m$^2$), $T^\mathrm{a}_t$ the air temperature, $\gamma = -0.004$~K$^{-1}$ the temperature coefficient of crystalline silicon~\citep{skoplaki2009temperature}, $T_\mathrm{NOCT} = 45$~\textdegree{}C (the median of the crystalline-silicon modules in the database of the California Energy Commission is $45.4$~\textdegree{}C~\citep{cec_modules}), and $L = 14\%$ the default system loss of PVGIS~\citep{pvgis_manual}.
The output is limited to the rated power, and hourly values are linearly interpolated to 15 minutes.
As a consistency check, the annual yield differs from that of the PVGIS hourly series for 2023~\citep{huld2012pvgis,pvgis_manual} by $-3.6\%$, $+1.7\%$ and $+1.5\%$ in Milan, Rome and Palermo, with hourly correlation $0.95$; this is a comparison between models, not with a measured array.
The PV forecast applies the same model to the irradiance and temperature forecasts of the ICON model~\citep{zangl2015icon} archived by Open-Meteo, which gives for each hour the value predicted $24$ hours earlier; successive hours can thus come from successive model runs, and \cref{sec:scenarios} states which values a policy may use at each decision.
The root-mean-square error of this forecast is $39$--$46\%$ lower than that of the persistence forecast (yesterday's profile).

\paragraph{Grid signals.}
Load and generation by production type in Italy are hourly data of the ENTSO-E Transparency Platform redistributed by Energy-Charts~\citep{energycharts}.
The average carbon intensity $\alpha_t$ is the generation-weighted average of operational emission factors ($900$~gCO$_2$/kWh for hard coal, $650$ for oil, $400$ for natural gas and coal-derived gas, $300$ for waste and other non-renewable sources, $0$ otherwise), excluding imports and the output of pumped storage.
These factors are consistent with the emission factors per unit of primary energy of Italian power plants~\citep[Table~1.2]{ispra2025emissioni} and typical efficiencies, and the resulting annual factor of gross production, $255.7$ and $214.8$~gCO$_2$/kWh in 2023 and 2024, matches the official values of $256.3$ and $215.9$~\citep[Table~1.3]{ispra2025emissioni}.
The marginal emission factor $\alpha^\mathrm{m}_t$ is estimated, following~\citet{hawkes2010marginal} and~\citet{silerevans2012marginal}, as the slope $a$ of the regression $\Delta E_h = a\, \Delta D_h + c$ of the hourly change of fossil emissions on the hourly change of the system, fitted separately for each year, quarter and hour of the day on the observations of that hour and of the two neighbouring ones (about $270$ per fit).
Unlike those studies, we regress on the residual load, i.e.\ load minus wind and solar generation, so that the variable renewables do not appear as changes of demand; resources that produce regardless of demand bias the estimate otherwise, as \citet{hawkes2010marginal} observed for run-of-river hydropower.
The median coefficient of determination is $0.60$ (between $0.3$ and $0.8$ for $80\%$ of the fits), and the estimates range between $49$ and $614$~gCO$_2$/kWh, so that the clipping to $[0, 1000]$ is never active.
The result is a table of $24$ values for each quarter, which is an observational estimate, not a causal measurement, and is estimated on the same quarter in which it is used; \cref{sec:closed-loop} also considers policies that plan with the table of the previous year.
The emissions attributed to the energy delivered by the battery are $\beta = 33$~gCO$_2$e/kWh, the median life-cycle value for lithium-ion storage of~\citet{nrel2021lca}; allocating them to throughput is an assumption, whose effect is examined in \cref{sec:headroom}.

\paragraph{Retail cost.}
The cost objective models an Italian resident customer with an offer indexed to the day-ahead price $\pi_t$ (\euro/MWh) of the Northern Italy bidding zone~\citep{energycharts}:
\begin{equation*}
    p^+_t = (1 + \nu) \left( (1 + \lambda) \pi_t / 1000 + c_\mathrm{v} \right), \qquad p^-_t = \pi_t / 1000,
\end{equation*}
where $\nu = 10\%$ is the VAT, $\lambda = 10\%$ the conventional loss factor of low-voltage withdrawals~\citep{arera584_2023}, and $c_\mathrm{v} = 0.0834$~\euro/kWh the sum of the per-kWh components for resident domestic customers in the first quarter of 2024: dispatching ($1.609$~c\euro/kWh~\citep{arera626_2023}), network ($1.057$~c\euro/kWh~\citep{arera630_2023}), system and equalization charges ($3.40$~c\euro/kWh~\citep{arera633_2023}) and excise duty ($2.27$~c\euro/kWh~\citep{adm2024accise}).
We keep these values for the whole year, and apply the excise duty to all imports, although resident supplies with up to $3$~kW are exempt for the first $150$~kWh of each month~\citep[art.~52]{dlgs504_1995}; fixed charges do not depend on the operation of the battery and are omitted.
Exported energy is paid at the hourly zonal price under the dedicated withdrawal scheme~\citep{arera280_2007}; we ignore the guaranteed minimum prices for small plants.
The wear cost is an assumption, $w = 0.05$~\euro/kWh delivered, with sensitivity $0$ and $0.10$~\euro/kWh; for reference, installed home storage costs $400$--$1{,}000$~\euro{} per kWh of usable capacity with a 15-year life~\citep{kost2024lcoe}.

\paragraph{Systems.}
We consider PV systems of $3$, $4.5$ and $6$~kWp and batteries of $5$~kWh / $2.5$~kW, $10$~kWh / $5$~kW and $15$~kWh / $5$~kW, with $c = 2$ cycles per day and constant efficiencies $\eta = \mu = \sqrt{0.96}$, i.e.\ an assumed round-trip efficiency of $96\%$ that neglects standby and part-load losses ($90\%$ in a sensitivity analysis).
The base case is $4.5$~kWp and $10$~kWh / $5$~kW in Milan.

\paragraph{Implementation.}
All models are written in Python~3.12 and solved with Gurobi~13.0~\citep{gurobi}, with barrier for PF and simplex for the lookahead models, which are updated and warm-started from one period to the next.
Each simulation runs as a single-threaded job on a computing cluster.

\section{Headroom of optimization}\label{sec:headroom}

\begin{figure}[t]
\centering
\includegraphics[width=\linewidth]{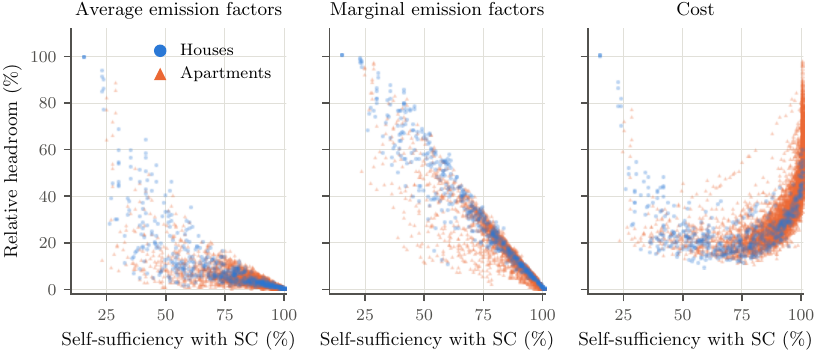}
\caption{Relative headroom $h$ against the self-sufficiency of the household with SC, for every household and PV/battery configuration in Milan ($9 \times 794$ points per panel). With emission factors the headroom vanishes as the rule covers the whole load; with cost it grows, because exports are paid and the rule leaves the battery idle.}
\label{fig:headroom}
\end{figure}

The gap between SC and PF bounds the benefit that any policy can obtain over the rule.
For each household we simulate SC and solve PF over the simulation horizon in $27$ site/PV/battery configurations for each emission objective and $9$ Milan configurations for the cost objective. The two wear-cost variants add $18$ cases, for a total of $81$ objective/configuration combinations per household and $64{,}314$ PF problems across the $794$ households.
Let $J_\mathrm{NB}$, $J_\mathrm{SC}$ and $J_\mathrm{PF}$ be the realized objectives over the $318$ days.
We report the relative headroom
\begin{equation*}
    h = \frac{J_\mathrm{SC} - J_\mathrm{PF}}{J_\mathrm{NB} - J_\mathrm{PF}},
\end{equation*}
the share of the benefit of a perfectly operated battery that the rule leaves on the table, and the absolute headroom $J_\mathrm{SC} - J_\mathrm{PF}$, in kgCO$_2$e or euros over the $318$ days.
The ratio $h$ is nonnegative and exceeds one when the rule is worse than no battery.
Medians are taken over households separately for each quantity, so that ``the median house'' below does not denote one physical household.
\Cref{tab:headroom} summarizes $h$ by objective, household type and configuration, and \cref{fig:headroom} shows it against the self-sufficiency of the household with SC, i.e.\ the share of the load not imported from the grid.

\begin{table}[t]
\caption{Relative headroom $h$ (\%) of optimization over the self-consumption rule: median over the households, with the interquartile range in the base configuration (Milan, $4.5$~kWp, $10$~kWh / $5$~kW), and median when one element of the base configuration is changed. Absolute headroom: median of $J_\mathrm{SC} - J_\mathrm{PF}$ over the $318$ days in the base configuration. The cost objective is computed for Milan only.}
\label{tab:headroom}
\centering
{\small
\setlength{\tabcolsep}{4pt}
\begin{tabular}{@{}lrrrrrrrr@{}}
\toprule
& \multicolumn{2}{c}{Base configuration} & \multicolumn{2}{c}{Site} & \multicolumn{2}{c}{PV} & \multicolumn{2}{c}{Battery} \\
\cmidrule(lr){2-3} \cmidrule(lr){4-5} \cmidrule(lr){6-7} \cmidrule(l){8-9}
Objective, households & $h$ (IQR) & Absolute & Rome & Palermo & 3 kWp & 6 kWp & 5 kWh & 15 kWh \\
\midrule
AEF, houses & 5.7 (2.7--12.6) & 19~kg & 2.0 & 2.2 & 17.3 & 3.2 & 5.5 & 6.7 \\
AEF, apartments & 0.6 (0.0--1.8) & 1~kg & 0.0 & 0.0 & 1.8 & 0.2 & 1.1 & 0.2 \\
\addlinespace
MEF, houses & 23.0 (11.2--42.1) & 137~kg & 14.3 & 11.9 & 43.3 & 16.8 & 33.5 & 22.1 \\
MEF, apartments & 2.1 (0.0--7.1) & 5~kg & 0.1 & 0.0 & 5.4 & 1.0 & 6.8 & 0.6 \\
\addlinespace
Cost, houses & 21.1 (16.6--28.8) & 36~\euro & -- & -- & 23.4 & 23.1 & 22.8 & 22.9 \\
Cost, apartments & 40.7 (30.3--53.5) & 41~\euro & -- & -- & 31.8 & 46.1 & 34.4 & 43.4 \\
\addlinespace
Cost, $w = 0$, houses & 26.6 (19.6--34.2) & 78~\euro & -- & -- & 32.7 & 24.9 & 22.5 & 29.6 \\
Cost, $w = 0$, apartments & 51.5 (38.4--65.1) & 110~\euro & -- & -- & 42.9 & 56.3 & 37.6 & 56.7 \\
\addlinespace
Cost, $w = 0.10$, houses & 38.7 (29.7--52.5) & 31~\euro & -- & -- & 38.6 & 46.7 & 43.4 & 39.9 \\
Cost, $w = 0.10$, apartments & 64.5 (52.4--78.2) & 31~\euro & -- & -- & 52.7 & 72.4 & 62.0 & 66.4 \\
\bottomrule
\end{tabular}
}
\end{table}

\paragraph{Emissions.}
With average emission factors the rule is close to optimal for most households.
In the base configuration the median relative headroom is $5.7\%$ for houses and $0.6\%$ for apartments, i.e.\ $19$ and $1$~kgCO$_2$e over the $318$ days, while the rule itself avoids $318$ and $176$~kg with respect to the household without battery.
With the estimated marginal factors the headroom of houses is four times larger, $23.0\%$ or $137$~kg, a third of the emissions of the median house with SC; that of apartments remains small ($2.1\%$).
The marginal factor varies much more within the day than the average one: over the simulation period, the median daily range is $410$~gCO$_2$/kWh against $115$.
A perfectly operated battery exploits this variation by charging from the grid when the factor is low: in the median house PF draws $450$~kWh from the grid into the battery with marginal factors and $301$~kWh with average ones.

\Cref{fig:headroom} shows that, for both emission objectives, the headroom is closely associated with the share of the load that the rule already covers.
Across all households and configurations, the rank correlation between $h$ and the self-sufficiency with SC is $-0.97$ with average and $-0.98$ with marginal factors (both quantities derive from the same import trajectories, so the correlation describes the data rather than testing a hypothesis).
When PV and battery cover almost the whole load there is little left to shift; the headroom grows with the energy that the household imports, which an optimized battery can buy at cleaner hours.
This explains the effect of the configuration in \cref{tab:headroom}: the headroom is largest with the smallest PV system and smaller in Rome and Palermo, where PV generation is higher, while the effect of the battery size is smaller and not consistent across objectives.
Over all configurations, $h$ is below $5\%$ for $94\%$ of the apartments and $56\%$ of the houses with average factors, and for $73\%$ and $22\%$ with marginal factors; it exceeds $20\%$ for $16\%$ of the houses with average factors, mostly with the smallest PV system.

\paragraph{Cost.}
With the cost objective the relative headroom is above $5\%$ for all households: the median is $21.1\%$ for houses and $40.7\%$ for apartments, or $36$ and $41$~\euro{} over the $318$ days, against $137$ and $61$~\euro{} saved by the rule with respect to no battery.
Unlike the emission objectives, the cost objective pays for exports, so a battery can also earn by selling stored PV energy when prices are high.
In the median apartment PF draws only $13$~kWh from the grid into the battery but cycles $0.47$ times a day, against $0.28$ with SC: most of the additional energy cycled is PV generation sold later in the day.
For this reason the relative headroom increases with self-sufficiency (rank correlation $0.88$): a household that rarely needs the battery at night leaves it idle under SC.
The wear cost changes the headroom considerably.
Without it, PF uses the battery for arbitrage, with $0.94$ cycles a day and $548$~kWh charged from the grid in the median house, and the absolute headroom more than doubles; with $w = 0.10$~\euro/kWh it falls to $31$~\euro{}, while the relative headroom increases because the benefit of the rule itself shrinks.
With this wear cost, the rule is even worse than no battery in some configurations of $22$ apartments with small loads and almost no imports, because the wear of the energy it cycles exceeds the value of shifting it (and, marginally, for one house with a very large load, whose battery is hardly used).
The objective includes the wear cost, so that these amounts are not savings on the electricity bill alone.

\paragraph{Sensitivity.}
We repeated the assessment in the base configuration under different assumptions on the battery.
With a round-trip efficiency of $90\%$ instead of $96\%$, the median relative headroom of houses is $5.3\%$, $24.3\%$ and $22.0\%$ with average factors, marginal factors and cost (against $5.7\%$, $23.0\%$ and $21.1\%$), and that of apartments changes by at most one point.
Without emissions attributed to the battery throughput ($\beta = 0$), or with twice the base value, the headroom of houses is $7.7\%$ and $4.5\%$ with average factors and $23.3\%$ and $23.8\%$ with marginal factors.
In the base configuration, the ranking of the headroom across objectives and household types is therefore retained under these efficiency and battery-emission assumptions.

They do depend on the accounting of exports.
If exports are credited at the marginal factor, so that the objective is the net exchange with the grid valued at $\alpha^\mathrm{m}_t$, the benefit of the rule over no battery falls to a median of $108$~kg for houses and $55$~kg for apartments in Milan ($136$ and $144$~kg for houses in Rome and Palermo), and the relative headroom rises to $86$--$90\%$ for houses and $95\%$ for apartments at the three sites.
Under this accounting PF operates the battery as an arbitrage device between hours of low and high marginal factor, with $1.6$ cycles a day and about $4{,}200$~kWh charged from the grid over the period, regardless of the household; in $4\%$ of the houses the rule is even worse than no battery.
This objective answers a different question, the effect on the emissions of the power system of a small battery that does not alter the marginal factors, and it rests even more heavily on the estimate of those factors; the closed-loop comparison includes it as a variant.

In summary, in these data an optimizing controller has little to gain with average emission factors for most households, a large potential with marginal factors in houses that rely on the grid, and a moderate but systematic potential with retail cost, with a median of a few tens of euros over the period.

\section{Closed-loop comparison}\label{sec:closed-loop}

\begin{figure}[t]
\centering
\includegraphics[width=\linewidth]{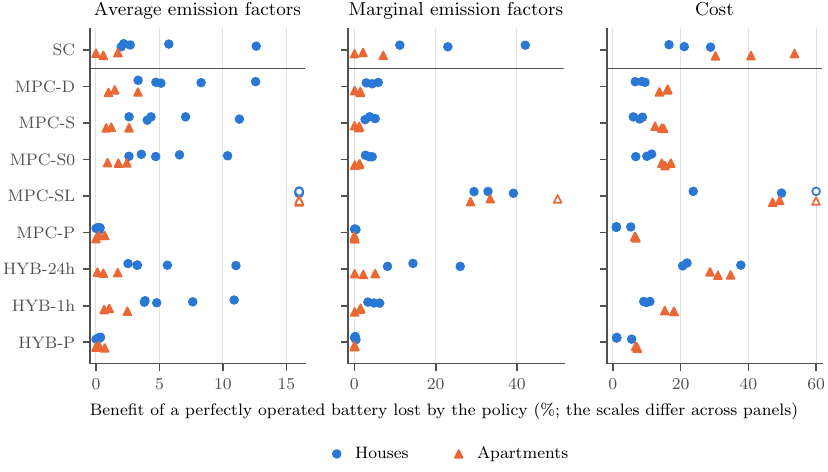}
\caption{Share of the benefit of a perfectly operated battery lost by each policy, $(J - J_\mathrm{PF}) / (J_\mathrm{NB} - J_\mathrm{PF})$, in the $20$ cases (one marker per case; for SC it is the relative headroom $h$). The horizontal scales differ across panels. Open markers at the right edge are beyond the axis: MPC-SL loses between $24\%$ and $86\%$.}
\label{fig:closed-loop}
\end{figure}

We select $20$ representative cases from \cref{sec:headroom}: in the base configuration, for each objective and separately for houses and apartments, the households closest to the 25th, 50th and 75th percentiles of the relative headroom; and the median house in Rome and Palermo with the average emission factor.
Each case is simulated with all the policies of \cref{tab:policies}.
Four variants are run on a subset of the cases.
For the three houses with marginal factors, the policies plan with the table of marginal factors estimated on the same quarter of the previous year and are evaluated with the actual one, so that no data of the evaluation year are used to fit the planning table; and the objective with exports credited at the marginal factor (\cref{sec:headroom}) is considered.
For the median house and the median apartment with the cost objective, the comparison is repeated with wear costs of $0$ and $0.10$~\euro/kWh.
\Cref{fig:closed-loop} shows, for each policy and case, the share $(J - J_\mathrm{PF}) / (J_\mathrm{NB} - J_\mathrm{PF})$ of the benefit of a perfectly operated battery that the policy loses, which equals $h$ for SC.
\Cref{tab:closed-loop} reports the share of the headroom recovered, $(J_\mathrm{SC} - J) / (J_\mathrm{SC} - J_\mathrm{PF})$, which is negative when the policy is worse than the rule.
Apartments with emission objectives are omitted from the table: their headroom, at most $3$~kg with average and $22$~kg with marginal factors, is too small for the ratio to be meaningful.
The groups contain three to five cases, and we give the range over the cases where the median alone would hide it.

\begin{table}[t]
\caption{Share of the SC--PF headroom recovered by each policy (\%), median over the cases of each group ($5$ houses with AEF, including Rome and Palermo; $3$ cases otherwise); negative values mean worse than SC. The first row is the median headroom $J_\mathrm{SC} - J_\mathrm{PF}$ over the $318$ days. The last columns give the number of linear programs solved per day and the median solver time per simulated day (s, one core).}
\label{tab:closed-loop}
\centering
{\small
\setlength{\tabcolsep}{5pt}
\begin{tabular}{@{}lrrrrrr@{}}
\toprule
& \multicolumn{2}{c}{Houses, emissions} & \multicolumn{2}{c}{Cost} & \multicolumn{2}{c}{Computation} \\
\cmidrule(lr){2-3} \cmidrule(lr){4-5} \cmidrule(l){6-7}
Policy & AEF & MEF & Houses & Apartments & LPs/day & s/day \\
\midrule
SC--PF headroom & 11~kg & 69~kg & 37~\euro & 41~\euro & & \\
\midrule
MPC-D & $-66$ & $81$ & $69$ & $66$ & 96 & 0.037 \\
MPC-S & $-31$ & $84$ & $71$ & $69$ & 96 & 9.7 \\
MPC-S0 & $-30$ & $84$ & $65$ & $65$ & 96 & 12 \\
MPC-SL & $<-100$ & $-70$ & $<-100$ & $-16$ & 96 & 10 \\
MPC-P & $95$ & $98$ & $94$ & $83$ & 96 & 0.036 \\
HYB-24h & $-21$ & $37$ & $-24$ & $24$ & 1 & 0.0007 \\
HYB-1h & $-44$ & $79$ & $53$ & $62$ & 24 & 0.017 \\
HYB-P & $95$ & $98$ & $93$ & $82$ & 24 & 0.016 \\
\bottomrule
\end{tabular}
}
\end{table}

\paragraph{Perfect information over the window.}
With the actual realizations over the next $24$ hours, MPC-P recovers $95$--$98\%$ of the headroom with emission objectives and $83$--$94\%$ with cost in the median of each group, and is better than SC in all $20$ cases.
The chosen 24-hour formulation can therefore recover most of the headroom when supplied with perfect information.
The remaining gap to PF can reflect the finite horizon, its terminal treatment and the resulting sequence of decisions.
Differences from forecast-based MPC also depend on the lookahead formulation and on the subsequent SOC trajectory; they do not isolate the forecast error alone.

\paragraph{Forecast-based MPC.}
With marginal factors, MPC-S recovers $76$--$88\%$ of the headroom in the three houses, i.e.\ $41$, $57$ and $334$~kg, and is better than SC in five of the six cases; in the apartment where it is worse, both policies are within $0.1$~kg of PF.
With cost it is better than SC in all six cases and recovers $47$--$73\%$ of the headroom, between $16$ and $32$~\euro{} over the period.
With average factors, instead, MPC-D, MPC-S, MPC-S0 and HYB-1h are all worse than SC in $7$ of the $8$ cases, the exception being the house with the largest headroom.
The headroom is small, $11$~kg in the median over the five houses, and forecast errors cost more.
In the median house in Milan, MPC-S charges from the grid about as much as PF ($434$ and $426$~kWh), but at the wrong times: the battery is then less able to store the PV surplus, and MPC-S imports $64$~kWh and exports $68$~kWh more than PF.
The losses are small in absolute terms, at most $8$~kg for MPC-S.
The two-stage model gains $2$--$3$ points over the deterministic one in the group medians with marginal factors and cost (MPC-S against MPC-D), at about $260$ times the solver time.

\paragraph{Value of forecasts and of measurements.}
Conditioning the scenarios on the weather forecast changes the median recovered share by $-1$ to $+6$ points (MPC-S against MPC-S0).
In these cases the forecast of the weather is thus worth little to the policy, in line with the observation of~\citet{litjens2018assessment}, for Dutch households and businesses, that a lower forecast error does not directly imply a higher storage revenue.
Measuring the current period, instead, is essential.
MPC-SL, which uses the same model and scenarios as MPC-S but decides before observing the period, loses $24$--$86\%$ of the benefit of a perfectly operated battery and is worse than SC in $17$ of the $20$ cases.
With delayed information the battery cannot adapt its decision to the realized net load of the current 15-minute period, and every mismatch is absorbed by the grid.
The simulation does not resolve fluctuations within a period, so that this result concerns information at the 15-minute scale only.

\paragraph{Hybrid policies.}
With a perfect plan, the setpoints lose almost nothing: HYB-P recovers within one point of MPC-P in every group.
With forecasts, HYB-1h is compared with MPC-D, which uses the same deterministic model every $15$ minutes.
With marginal factors it recovers $70$--$85\%$ of the headroom in the three houses, against $74$--$86\%$ for MPC-D.
With cost it recovers $34$--$68\%$ in the six cases, against $43$--$70\%$ for MPC-D; the difference is largest in the median house ($53\%$ against $69\%$).
It is worse than SC in the same $8$ cases as MPC-D and MPC-S.
It solves $24$ small linear programs a day instead of $96$, with a median solver time of $0.017$~s per simulated day against $0.037$~s for MPC-D and $9.7$~s for MPC-S.
On one core of a cluster node, a full simulation of $318$ days took a median of $22$~s with HYB-1h, $58$~s with MPC-D and $56$~min with MPC-S, including the generation of the scenarios.
All these times are negligible with respect to the 15-minute period, so that the difference matters for the hardware and software that a controller requires rather than for feasibility.
A plan computed once a day is not enough: HYB-24h recovers $37\%$ in the median house with marginal factors and is worse than SC in $10$ of the $20$ cases.

\begin{figure}[t]
\centering
\includegraphics[width=\linewidth]{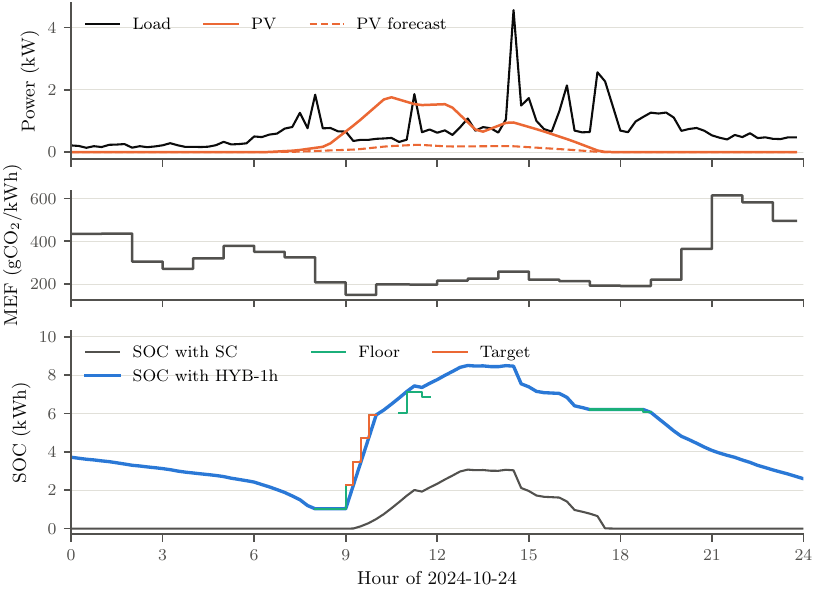}
\caption{Operation of HYB-1h on 24 October 2024 in the median house with marginal emission factors (Milan, $4.5$~kWp, $10$~kWh / $5$~kW). Top: load, PV generation and its day-ahead forecast. Middle: marginal emission factor. Bottom: SOC with SC and with HYB-1h, and the setpoints of the hybrid policy (shown only when active); the two policies start the day from different SOC because of their different history.}
\label{fig:day}
\end{figure}

\Cref{fig:day} illustrates HYB-1h on a day of October in the median house with marginal factors.
The forecast underestimated the PV generation considerably.
The plan sets a target in the morning, when the marginal factor is low, and the battery charges from the grid; it then sets a floor that holds energy for the evening, when the factor reaches $600$~gCO$_2$/kWh.
Between these setpoints the rule stores the PV surplus that the forecast did not anticipate, as SC would.

\begin{table}[t]
\caption{Variants of the closed-loop comparison: headroom $J_\mathrm{SC} - J_\mathrm{PF}$ over the $318$ days and share recovered by four policies (\%). Marginal factors: range over the three houses. Cost: median house and median apartment with three wear costs $w$ (\euro/kWh; $0.05$ is the base value). MPC-P was not run with the previous-year factors.}
\label{tab:variants}
\centering
{\small
\setlength{\tabcolsep}{5pt}
\begin{tabular}{@{}lrrrrr@{}}
\toprule
Variant & Headroom & MPC-D & MPC-S & MPC-P & HYB-1h \\
\midrule
MEF, base & 54--380~kg & $74$--$86$ & $76$--$88$ & $98$--$100$ & $70$--$85$ \\
MEF, previous-year factors & 54--380~kg & $65$--$73$ & $70$--$74$ & -- & $66$--$71$ \\
MEF, exports credited & 880--1000~kg & $100$ & $100$ & $100$ & $99$ \\
\addlinespace
Cost, $w = 0$, house & 97~\euro & $73$ & $78$ & $100$ & $75$ \\
Cost, $w = 0.05$, house & 43~\euro & $69$ & $71$ & $95$ & $53$ \\
Cost, $w = 0.10$, house & 37~\euro & $76$ & $77$ & $92$ & $46$ \\
\addlinespace
Cost, $w = 0$, apartment & 109~\euro & $89$ & $89$ & $98$ & $87$ \\
Cost, $w = 0.05$, apartment & 41~\euro & $66$ & $69$ & $83$ & $62$ \\
Cost, $w = 0.10$, apartment & 32~\euro & $58$ & $59$ & $71$ & $51$ \\
\bottomrule
\end{tabular}
}
\end{table}

\paragraph{Variants.}
\Cref{tab:variants} reports the variants.
In the base runs the policies know the table of marginal factors estimated on the evaluation period itself.
When they plan with the table of the same quarter of the previous year, MPC-S recovers $70$--$74\%$ of the headroom in the three houses instead of $76$--$88\%$, a loss of $6$ to $18$ points, and remains better than SC in all three; MPC-D and HYB-1h behave similarly.
Most of the benefit with marginal factors therefore does not depend on knowing the factors of the evaluation year.

When exports are credited at the marginal factor, MPC-D and MPC-S attain the perfect-foresight value in the three tested cases, and HYB-1h recovers $99\%$ of the headroom.
The reason lies in the objective rather than in the policies: the net exchange valued at $\alpha^\mathrm{m}_t$ is the sum of a term that does not depend on the decisions and of the battery flows valued at $\alpha^\mathrm{m}_t$, so that the problem reduces to arbitrage on a signal that the policies know, and the forecasts of load and PV generation play no role in the battery decisions; a finite lookahead could still lose against PF through its terminal effects, which here proved negligible.
Under this accounting the difficulty moves entirely to the forecast of the marginal factors, which this variant does not test.

Without wear cost the headroom more than doubles, and MPC-S recovers $78\%$ and $89\%$ of it in the median house and apartment; with $w = 0.10$~\euro/kWh it recovers $77\%$ and $59\%$ of a smaller headroom.
HYB-1h is close to MPC-D without wear cost, and falls behind as the wear cost grows, in the house from $75\%$ against $73\%$ to $46\%$ against $76\%$.
A likely reason is that the setpoints let the rule cycle energy that the plan, which accounts for wear, would not cycle.

\section{Discussion}\label{sec:discussion}

\paragraph{When the rule is enough.}
With average emission factors and import-based accounting, the rule is hard to beat in these data: the median headroom is a few percent of the benefit of the battery, and the forecast-based policies we tested were worse than the rule in $7$ of the $8$ simulated cases.
The same holds for most apartments with either emission factor, whose load the rule covers almost entirely.
The headroom of the emission objectives is closely associated with the self-sufficiency obtained with the rule (\cref{fig:headroom}), which can be computed from the meter data of an existing installation, or from a simulation of the rule, without solving any optimization problem.
It is therefore a simple first screen of the households for which a more complex controller can matter.

\paragraph{When optimization can help.}
With marginal emission factors, houses that import a significant share of their load have a large headroom: in the base configuration the median is $137$~kg, a third of the emissions of the median house with the rule.
In the three houses simulated in closed loop, scenario-based MPC recovers $76$--$88\%$ of the headroom, and $70$--$74\%$ when it plans with the factors of the previous year.
With cost, all households of the data set have a headroom, with a median of a few tens of euros over the period, of which forecast-based MPC recovers between a half and three quarters in the six simulated cases.
These amounts are small compared with the investment in a battery, which costs $400$--$1{,}000$~\euro{} per kWh of usable capacity~\citep{kost2024lcoe}, and we do not assess the cost of the controller; they are also sensitive to the wear cost, which determines how much arbitrage is worth.
The emission results hold for the emissions of imported electricity: if exports are credited at the marginal factor, the battery becomes an arbitrage device whose operation is independent of the household, and the question becomes one of forecasting the marginal factors.
Which signal to use is a choice of the operator: as the literature on cost-driven operation shows~\citep{babacan2018unintended,olivieri2020optimization}, a battery operated for cost does not necessarily reduce emissions, and a battery operated for average factors gains little over the rule.

\paragraph{How to implement it.}
Measuring the current period was worth more than the weather forecast in our cases.
A policy that fixes the battery power before the period starts was worse than the plain rule in $17$ of the $20$ cases, even with the same model and scenarios as the best MPC, so the final decision should rest with a controller that acts on the measured net load.
The hybrid policy is one such controller.
With marginal factors it recovered nearly as much as deterministic MPC with the same forecasts ($70$--$85\%$ against $74$--$86\%$ of the headroom in the three houses), solving one small linear program per hour and keeping as real-time layer the rule already running in residential systems.
With cost it recovered less ($34$--$68\%$ against $43$--$70\%$), and the difference grew with the wear cost.
Because the plan reaches the battery only as four numbers per period, it could be computed remotely, and the controller could revert to neutral setpoints, i.e.\ to the plain rule, if the communication fails; we did not simulate such failures.
The plan has to be updated during the day, however: a daily plan was worse than the rule in half of the cases.

\paragraph{Limitations.}
The case study combines sources: the loads are measured in Switzerland, in households without electric heating or vehicle charging, and are paired with Italian weather and grid data of the same days.
The calendar is aligned, but the response of Italian households to the local weather, in particular for cooling, and Italian holidays are not represented.
The sample has $49$ houses and $745$ apartments of one utility, and the base battery is large for the median apartment; the results are a conditional comparison of operating policies, not a statement on which systems households should install.
The PV generation is computed from satellite irradiance rather than measured, and hourly irradiance interpolated to 15 minutes is smoother than the output of a real array; real-time information is idealized at the 15-minute resolution of the data.
The emission objectives count imports and do not credit exports (\cref{sec:headroom} shows how much this matters), the average factors exclude imported electricity, and the marginal factors are an observational estimate.
Battery degradation is represented only by the cycle limit, a constant efficiency and, for the cost objective, an assumed wear cost; limits on the power exchanged with the grid are not modelled.
The closed-loop comparison covers $20$ cases, selected by their headroom in the base configuration; they span that criterion but are not a random sample of the households, and the medians over three to five cases are not population estimates.
The average emission factor is forecast by our scenario model from its recent history; a dedicated forecast of the carbon intensity could change the results of the forecast-based policies with this objective, as could a different lookahead horizon or terminal condition.
The study covers $318$ days and a single household at a time; aggregation and grid services are not considered.
All results are simulations: unlike real-world validations of home energy management systems~\citep{ruddick2024realworld}, the policies have not been tested on hardware.

\section{Conclusions}\label{sec:conclusions}

We asked when replacing the self-consumption rule of a residential PV-battery system with an optimization-based policy improves emissions or cost, using measured loads, satellite-based PV generation with archived weather forecasts, and the Italian grid.
The answer depends on the objective and on how much of the load the rule already covers.
For the emissions of imported electricity valued at average factors, the rule is near-optimal for most households, and the forecast-based policies we tested did worse than it in $7$ of $8$ cases.
With marginal emission factors the headroom is large for houses that rely on the grid, and with retail cost its median is a few tens of euros over the period; in the simulated cases, scenario-based MPC recovered $76$--$88\%$ of the headroom of the houses with marginal factors and $47$--$73\%$ with cost.
Measuring the current period was essential, while the weather forecast added little.
A hybrid policy, in which the rule executes four SOC setpoints from an hourly linear-programming plan, recovered $70$--$85\%$ of the headroom with marginal factors and $34$--$68\%$ with cost, with one small linear program per hour and the rule as a fallback.
The conclusions on emissions depend on the accounting of exports, and all of them on a case study that pairs Swiss loads with Italian weather and grid data.
Natural extensions are dedicated forecasts of the grid signals, in particular of marginal factors, the operation of many batteries together, and the validation of the hybrid policy on a real installation.

\appendix

\section{Scenario model}\label{app:scenarios}

Let $v_d \in \reals^{n_v}$, $n_v = JT$, stack, for training day $d$, the $J = 4$ profiles of load, PV generation, $\log \max(x_t, 1)$ and PV forecast, where $x_t$ is the grid signal (gCO$_2$/kWh or \euro/MWh); MPC-S0 omits the last block ($J = 3$).
Let $\bar{v}$ be the average over the $n = 28$ training days and $D$ the diagonal matrix that holds, for each block, the standard deviation of all its centred entries.
The standardized vectors $u_d = D^{-1} (v_d - \bar{v})$ have sample covariance with eigenvalues $\lambda_1 \ge \lambda_2 \ge \cdots$ (at most $n - 1$ are nonzero) and eigenvectors $q_i$.
Following~\citet{tipping1999probabilistic}, we set $u = W z + \varepsilon$ with $z \sim \mathcal{N}(0, I_K)$, $\varepsilon \sim \mathcal{N}(0, \sigma^2 I)$,
\begin{equation*}
    \sigma^2 = \frac{1}{n_v - K} \sum_{i > K} \lambda_i, \qquad W = \left[ \sqrt{\lambda_1 - \sigma^2}\, q_1, \ldots, \sqrt{\lambda_K - \sigma^2}\, q_K \right],
\end{equation*}
and $K = 5$.
For numerical safety, $10^{-9}$ is added to the block standard deviations in $D$, and $\sigma^2$ and $\lambda_i - \sigma^2$ are bounded below by $10^{-8}$ and $10^{-12}$.
Given the observed entries $u_o$ of a day, with $W_o$ the corresponding rows of $W$, the latent coefficients have the Gaussian posterior
\begin{equation*}
    z \mid u_o \sim \mathcal{N}\left( \Sigma W_o^T u_o / \sigma^2,\ \Sigma \right), \qquad \Sigma = \left( I_K + W_o^T W_o / \sigma^2 \right)^{-1}.
\end{equation*}
A scenario is $\bar{v} + D W z$ with $z$ sampled from the posterior, mapped back to the original units.
The residual noise $\varepsilon$ is not added to the reconstruction, so that the scenarios are smooth profiles in the span of the principal components and understate the variability at the shortest time scales.
Load and PV generation are truncated at zero, PV generation is limited to the rated power and set to zero in the periods in which it is zero on all training days, and the observed entries are returned as observed.
The window of a lookahead model is composed of the remaining periods of the current day and the first periods of the next one, whose latent coefficients are sampled independently, each conditioned on its own observed entries (\cref{sec:scenarios}); no dependence across midnight is modelled beyond the common training window.

\section*{Declaration of competing interest}
The authors declare that they have no known competing financial interests or personal relationships that could have appeared to influence the work reported in this paper.

\section*{Funding}
This research did not receive any specific grant from funding agencies in the public, commercial, or not-for-profit sectors.

\section*{Data availability}
All data used in this study are public.
The household loads are from~\citet{kaiser2026swissloads} (CC BY 4.0); the Italian generation, load and day-ahead prices from Energy-Charts~\citep{energycharts}; the satellite irradiance, reanalysis temperature and archived weather forecasts from Open-Meteo~\citep{zippenfenig2023openmeteo} (CC BY 4.0).
The code and the processed data are available at~\citet{bernardelligithub}, except the processed household loads, which the code rebuilds from the public data set.

\section*{Declaration of generative AI and AI-assisted technologies in the manuscript preparation process}
During the preparation of this work, the authors used Claude (Anthropic) to assist with writing code for data processing, conducting experiments, and generating LaTeX scripts for figures.
After using this tool, the authors reviewed and edited the content as needed and take full responsibility for the content of the publication.

\bibliographystyle{plainnat}
\bibliography{references}

\begin{thebibliography}{51}
\providecommand{\natexlab}[1]{#1}
\providecommand{\url}[1]{\texttt{#1}}
\expandafter\ifx\csname urlstyle\endcsname\relax
  \providecommand{\doi}[1]{doi: #1}\else
  \providecommand{\doi}{doi: \begingroup \urlstyle{rm}\Url}\fi

\bibitem[{Agenzia delle Dogane e dei Monopoli}(2024)]{adm2024accise}
{Agenzia delle Dogane e dei Monopoli}.
\newblock Aliquote di imposta vigenti nel settore delle accise: aggiornamento
  al 1 gennaio 2024.
\newblock
  \url{https://adm.gov.it/portale/documents/20182/43975520/Aliquote-ACCISE-naz-01012024.pdf},
  2024.

\bibitem[Aikman et~al.(2026)Aikman, Munz, Tabares-Velasco, and
  Newman]{aikman2026sizing}
Andrew Aikman, Karlyle Munz, Paulo~Cesar Tabares-Velasco, and Alexandra~M.
  Newman.
\newblock Sizing of stand-alone residential battery energy storage systems for
  energy cost and emissions savings.
\newblock \emph{Journal of Energy Storage}, 148:\penalty0 119541, 2026.
\newblock \doi{10.1016/j.est.2025.119541}.

\bibitem[Amaripadath et~al.(2026)Amaripadath, Sailor, Khan, Bertini, and
  Attia]{amaripadath2026achieving}
Deepak Amaripadath, David~J. Sailor, Ansar Khan, Aurora Bertini, and Shady
  Attia.
\newblock Achieving low-carbon electricity operations in net zero energy
  residential buildings through plug-in battery energy storage systems.
\newblock \emph{Journal of Energy Storage}, 153:\penalty0 120964, 2026.
\newblock \doi{10.1016/j.est.2026.120964}.

\bibitem[Angenendt et~al.(2018)Angenendt, Zurmühlen, Axelsen, and
  Sauer]{angenendt2018comparison}
Georg Angenendt, Sebastian Zurmühlen, Hendrik Axelsen, and Dirk~Uwe Sauer.
\newblock Comparison of different operation strategies for pv battery home
  storage systems including forecast-based operation strategies.
\newblock \emph{Applied Energy}, 229:\penalty0 884–899, 2018.
\newblock \doi{10.1016/j.apenergy.2018.08.058}.

\bibitem[{ARERA}(2007)]{arera280_2007}
{ARERA}.
\newblock Deliberazione 6 novembre 2007, n. 280/07: Modalità e condizioni
  tecnico-economiche per il ritiro dell'energia elettrica ai sensi
  dell'articolo 13, commi 3 e 4, del decreto legislativo 29 dicembre 2003, n.
  387/03, e del comma 41 della legge 23 agosto 2004, n. 239/04.
\newblock \url{https://www.arera.it/atti-e-provvedimenti/dettaglio/07/280-07},
  2007.

\bibitem[{ARERA}(2023{\natexlab{a}})]{arera584_2023}
{ARERA}.
\newblock Deliberazione 12 dicembre 2023, 584/2023/{R/EEL}: Regolazione delle
  perdite sulle reti di distribuzione di energia elettrica.
\newblock \url{https://www.arera.it/fileadmin/allegati/docs/23/584-24.pdf},
  2023{\natexlab{a}}.

\bibitem[{ARERA}(2023{\natexlab{b}})]{arera626_2023}
{ARERA}.
\newblock Deliberazione 28 dicembre 2023, 626/2023/{R/EEL}: Aggiornamento, per
  il trimestre 1 gennaio -- 31 marzo 2024, delle condizioni economiche del
  servizio di vendita dell'energia elettrica in maggior tutela.
\newblock \url{https://www.arera.it/atti-e-provvedimenti/dettaglio/23/626-23},
  2023{\natexlab{b}}.

\bibitem[{ARERA}(2023{\natexlab{c}})]{arera630_2023}
{ARERA}.
\newblock Deliberazione 28 dicembre 2023, 630/2023/{R/EEL}: Aggiornamento, per
  l'anno 2024, delle tariffe per l'uso delle infrastrutture per i servizi di
  distribuzione e misura dell'energia elettrica per i clienti domestici, non
  domestici e delle condizioni economiche per l'erogazione del servizio di
  connessione.
\newblock \url{https://www.arera.it/atti-e-provvedimenti/dettaglio/23/630-23},
  2023{\natexlab{c}}.

\bibitem[{ARERA}(2023{\natexlab{d}})]{arera633_2023}
{ARERA}.
\newblock Deliberazione 28 dicembre 2023, 633/2023/{R/COM}: Aggiornamento, dal
  1 gennaio 2024, delle componenti tariffarie destinate alla copertura degli
  oneri generali e di ulteriori componenti del settore elettrico e del settore
  gas.
\newblock \url{https://www.arera.it/atti-e-provvedimenti/dettaglio/23/633-23},
  2023{\natexlab{d}}.

\bibitem[{ARERA}(2024)]{arera2024schedaQ3}
{ARERA}.
\newblock Scheda tecnica: l'aggiornamento delle condizioni di tutela, {III}
  trimestre 2024 nel dettaglio.
\newblock \url{https://www.arera.it/fileadmin/allegati/schede/240627st.pdf},
  2024.
\newblock 27 June 2024.

\bibitem[Azuatalam et~al.(2019)Azuatalam, Paridari, Ma, Förstl, Chapman, and
  Verbič]{azuatalam2019energy}
Donald Azuatalam, Kaveh Paridari, Yiju Ma, Markus Förstl, Archie~C. Chapman,
  and Gregor Verbič.
\newblock Energy management of small-scale pv-battery systems: A systematic
  review considering practical implementation, computational requirements,
  quality of input data and battery degradation.
\newblock \emph{Renewable and Sustainable Energy Reviews}, 112:\penalty0
  555–570, 2019.
\newblock \doi{10.1016/j.rser.2019.06.007}.

\bibitem[Babacan et~al.(2018)Babacan, Abdulla, Hanna, Kleissl, and
  Victor]{babacan2018unintended}
Oytun Babacan, Ahmed Abdulla, Ryan Hanna, Jan Kleissl, and David~G. Victor.
\newblock Unintended effects of residential energy storage on emissions from
  the electric power system.
\newblock \emph{Environmental Science \& Technology}, 52\penalty0
  (22):\penalty0 13600--13608, 2018.

\bibitem[Beaudin and Zareipour(2015)]{beaudin2015home}
Marc Beaudin and Hamidreza Zareipour.
\newblock Home energy management systems: A review of modelling and complexity.
\newblock \emph{Renewable and Sustainable Energy Reviews}, 45:\penalty0
  318–335, 2015.
\newblock \doi{10.1016/j.rser.2015.01.046}.

\bibitem[Bernardelli(2026)]{bernardelligithub}
Ambrogio~Maria Bernardelli.
\newblock Linear battery models.
\newblock GitHub repository, version bab4603, 2026.
\newblock
  \url{https://github.com/AmbrogioMB/linear-battery-models/tree/bab4603653d1fb434ba61a31c6f7fe4f879fe0c1}.

\bibitem[{California Energy Commission}(2026)]{cec_modules}
{California Energy Commission}.
\newblock Solar equipment lists: photovoltaic modules ({CEC} module database,
  as distributed in the {NREL} system advisor model library ``cec
  modules.csv'', version 2026.7.3).
\newblock \url{https://github.com/NREL/SAM/tree/develop/deploy/libraries},
  2026.
\newblock Accessed 13 November 2026.

\bibitem[Caputo(2025)]{ispra2025emissioni}
Antonio Caputo.
\newblock Le emissioni di {CO2} nel settore elettrico nazionale e regionale.
\newblock Rapporti 413/2025, ISPRA -- Istituto Superiore per la Protezione e la
  Ricerca Ambientale, 2025.

\bibitem[Elkazaz et~al.(2020)Elkazaz, Sumner, Naghiyev, Pholboon, Davies, and
  Thomas]{elkazaz2020hierarchical}
Mahmoud Elkazaz, Mark Sumner, Eldar Naghiyev, Seksak Pholboon, Richard Davies,
  and David Thomas.
\newblock A hierarchical two-stage energy management for a home microgrid using
  model predictive and real-time controllers.
\newblock \emph{Applied Energy}, 269:\penalty0 115118, 2020.
\newblock \doi{10.1016/j.apenergy.2020.115118}.

\bibitem[{European Commission, Joint Research Centre}(n.d.)]{pvgis_manual}
{European Commission, Joint Research Centre}.
\newblock {PVGIS} user manual.
\newblock
  \url{https://joint-research-centre.ec.europa.eu/photovoltaic-geographical-information-system-pvgis/getting-started-pvgis/pvgis-user-manual_en},
  n.d.
\newblock Accessed 13 November 2026.

\bibitem[Fisher and Apt(2017)]{fisher2017emissions}
Michael~J. Fisher and Jay Apt.
\newblock Emissions and economics of behind-the-meter electricity storage.
\newblock \emph{Environmental Science \& Technology}, 51\penalty0 (3):\penalty0
  1094–1101, 2017.
\newblock \doi{10.1021/acs.est.6b03536}.

\bibitem[{Fraunhofer Institute for Solar Energy Systems
  ISE}(2026)]{energycharts}
{Fraunhofer Institute for Solar Energy Systems ISE}.
\newblock {Energy-Charts}: Public net electricity generation and load in
  {I}taly, 2026.
\newblock URL \url{https://www.energy-charts.info}.

\bibitem[{Gurobi Optimization, LLC}(2026)]{gurobi}
{Gurobi Optimization, LLC}.
\newblock {Gurobi Optimizer Reference Manual}, version 13.0, 2026.
\newblock URL \url{https://www.gurobi.com}.

\bibitem[Hafiz et~al.(2020)Hafiz, Awal, Queiroz, and Husain]{hafiz2020realtime}
Faeza Hafiz, M.~A. Awal, Anderson Rodrigo~de Queiroz, and Iqbal Husain.
\newblock Real-time stochastic optimization of energy storage management using
  deep learning-based forecasts for residential pv applications.
\newblock \emph{IEEE Transactions on Industry Applications}, 56\penalty0
  (3):\penalty0 2216–2226, 2020.
\newblock \doi{10.1109/tia.2020.2968534}.

\bibitem[Hanna et~al.(2014)Hanna, Kleissl, Nottrott, and
  Ferry]{hanna2014energy}
R.~Hanna, J.~Kleissl, A.~Nottrott, and M.~Ferry.
\newblock Energy dispatch schedule optimization for demand charge reduction
  using a photovoltaic-battery storage system with solar forecasting.
\newblock \emph{Solar Energy}, 103:\penalty0 269–287, 2014.
\newblock \doi{10.1016/j.solener.2014.02.020}.

\bibitem[Hawkes(2010)]{hawkes2010marginal}
A.~D. Hawkes.
\newblock Estimating marginal {CO2} emissions rates for national electricity
  systems.
\newblock \emph{Energy Policy}, 38:\penalty0 5977--5987, 2010.
\newblock \doi{10.1016/j.enpol.2010.05.053}.

\bibitem[Hersbach et~al.(2020)Hersbach, Bell, Berrisford, Hirahara,
  et~al.]{hersbach2020era5}
Hans Hersbach, Bill Bell, Paul Berrisford, Shoji Hirahara, et~al.
\newblock The {ERA5} global reanalysis.
\newblock \emph{Quarterly Journal of the Royal Meteorological Society},
  146:\penalty0 1999--2049, 2020.
\newblock \doi{10.1002/qj.3803}.

\bibitem[Huld et~al.(2012)Huld, M{\"u}ller, and Gambardella]{huld2012pvgis}
Thomas Huld, Richard M{\"u}ller, and Attilio Gambardella.
\newblock A new solar radiation database for estimating {PV} performance in
  {Europe} and {Africa}.
\newblock \emph{Solar Energy}, 86:\penalty0 1803--1815, 2012.
\newblock \doi{10.1016/j.solener.2012.03.006}.

\bibitem[Kaiser et~al.(2026)Kaiser, Kreft, Gonz{\'a}lez~Vay{\'a}, and
  Hug]{kaiser2026swissloads}
Katharina Kaiser, Markus Kreft, Marina Gonz{\'a}lez~Vay{\'a}, and Gabriela Hug.
\newblock Dataset on residential electricity load profiles in {Switzerland}.
\newblock \url{https://doi.org/10.5281/zenodo.21398768}, 2026.

\bibitem[Kazhamiaka et~al.(2016)Kazhamiaka, Rosenberg, and
  Keshav]{kazhamiaka2016practical}
Fiodar Kazhamiaka, Catherine Rosenberg, and Srinivasan Keshav.
\newblock Practical strategies for storage operation in energy systems: Design
  and evaluation.
\newblock \emph{IEEE Transactions on Sustainable Energy}, 7\penalty0
  (4):\penalty0 1602–1610, 2016.
\newblock \doi{10.1109/tste.2016.2569425}.

\bibitem[Kost et~al.(2024)Kost, M{\"u}ller, Sep{\'u}lveda~Schweiger, Fluri, and
  Thomsen]{kost2024lcoe}
Christoph Kost, Paul M{\"u}ller, Jael Sep{\'u}lveda~Schweiger, Verena Fluri,
  and Jessica Thomsen.
\newblock Stromgestehungskosten erneuerbare energien.
\newblock Technical report, Fraunhofer-Institut f{\"u}r Solare Energiesysteme
  ISE, July 2024.

\bibitem[Langenmayr et~al.(2020)Langenmayr, Wang, and
  Jochem]{langenmayr2020unit}
Uwe Langenmayr, Weimin Wang, and Patrick Jochem.
\newblock Unit commitment of photovoltaic-battery systems: An advanced approach
  considering uncertainties from load, electric vehicles, and photovoltaic.
\newblock \emph{Applied Energy}, 280:\penalty0 115972, 2020.
\newblock \doi{10.1016/j.apenergy.2020.115972}.

\bibitem[Litjens et~al.(2018)Litjens, Worrell, and van
  Sark]{litjens2018assessment}
G.B.M.A. Litjens, E.~Worrell, and W.G.J.H.M. van Sark.
\newblock Assessment of forecasting methods on performance of
  photovoltaic-battery systems.
\newblock \emph{Applied Energy}, 221:\penalty0 358–373, 2018.
\newblock \doi{10.1016/j.apenergy.2018.03.154}.

\bibitem[Luthander et~al.(2015)Luthander, Wid{\'e}n, Nilsson, and
  Palm]{luthander2015photovoltaic}
Rasmus Luthander, Joakim Wid{\'e}n, Daniel Nilsson, and Jenny Palm.
\newblock Photovoltaic self-consumption in buildings: A review.
\newblock \emph{Applied Energy}, 142:\penalty0 80--94, 2015.

\bibitem[Mesbah(2016)]{mesbah2016stochastic}
Ali Mesbah.
\newblock Stochastic model predictive control: An overview and perspectives for
  future research.
\newblock \emph{IEEE Control Systems Magazine}, 36\penalty0 (6):\penalty0
  30--44, 2016.

\bibitem[{National Renewable Energy Laboratory}(2021)]{nrel2021lca}
{National Renewable Energy Laboratory}.
\newblock Life cycle greenhouse gas emissions from electricity generation:
  Update.
\newblock Technical Report NREL/FS-6A50-80580, NREL, 2021.

\bibitem[Olivieri and McConky(2020)]{olivieri2020optimization}
Zachary~T. Olivieri and Katie McConky.
\newblock Optimization of residential battery energy storage system scheduling
  for cost and emissions reductions.
\newblock \emph{Energy and Buildings}, 210:\penalty0 109787, 2020.
\newblock \doi{10.1016/j.enbuild.2020.109787}.

\bibitem[Parisio et~al.(2014)Parisio, Rikos, and Glielmo]{parisio2014model}
Alessandra Parisio, Evangelos Rikos, and Luigi Glielmo.
\newblock A model predictive control approach to microgrid operation
  optimization.
\newblock \emph{IEEE Transactions on Control Systems Technology}, 22\penalty0
  (5):\penalty0 1813--1827, 2014.

\bibitem[Pfeifroth et~al.(2023)Pfeifroth, Kothe, Dr{\"u}cke, Trentmann,
  et~al.]{pfeifroth2023sarah3}
Uwe Pfeifroth, Steffen Kothe, Jaqueline Dr{\"u}cke, J{\"o}rg Trentmann, et~al.
\newblock Surface radiation data set -- heliosat ({SARAH}) -- edition 3.
\newblock \url{https://doi.org/10.5676/EUM_SAF_CM/SARAH/V003}, 2023.

\bibitem[Powell(2019)]{powell2019unified}
Warren~B. Powell.
\newblock A unified framework for stochastic optimization.
\newblock \emph{European Journal of Operational Research}, 275\penalty0
  (3):\penalty0 795--821, 2019.

\bibitem[Powell and Meisel(2016)]{powell2016tutorial}
Warren~B. Powell and Stephan Meisel.
\newblock Tutorial on stochastic optimization in energy---{P}art {I}: Modeling
  and policies.
\newblock \emph{IEEE Transactions on Power Systems}, 31\penalty0 (2):\penalty0
  1459--1467, 2016.

\bibitem[Rawlings et~al.(2017)Rawlings, Mayne, and Diehl]{rawlings2017model}
James~B. Rawlings, David~Q. Mayne, and Moritz~M. Diehl.
\newblock \emph{Model Predictive Control: Theory, Computation, and Design}.
\newblock Nob Hill Publishing, Madison, WI, 2nd edition, 2017.

\bibitem[{Repubblica Italiana}(1995)]{dlgs504_1995}
{Repubblica Italiana}.
\newblock Decreto legislativo 26 ottobre 1995, n. 504: Testo unico delle
  disposizioni legislative concernenti le imposte sulla produzione e sui
  consumi e relative sanzioni penali e amministrative.
\newblock
  \url{https://www.normattiva.it/uri-res/N2Ls?urn:nir:stato:decreto.legislativo:1995-10-26;504},
  1995.

\bibitem[Riffonneau et~al.(2011)Riffonneau, Bacha, Barruel, and
  Ploix]{riffonneau2011optimal}
Y~Riffonneau, S~Bacha, F~Barruel, and S~Ploix.
\newblock Optimal power flow management for grid connected pv systems with
  batteries.
\newblock \emph{IEEE Transactions on Sustainable Energy}, 2\penalty0
  (3):\penalty0 309–320, 2011.
\newblock \doi{10.1109/tste.2011.2114901}.

\bibitem[Ruddick et~al.(2024)Ruddick, Ceusters, Van~Kriekinge, Genov,
  De~Cauwer, Coosemans, and Messagie]{ruddick2024realworld}
Julian Ruddick, Glenn Ceusters, Gilles Van~Kriekinge, Evgenii Genov, Cedric
  De~Cauwer, Thierry Coosemans, and Maarten Messagie.
\newblock Real-world validation of safe reinforcement learning, model
  predictive control and decision tree-based home energy management systems.
\newblock \emph{Energy and AI}, 18:\penalty0 100448, 2024.
\newblock \doi{10.1016/j.egyai.2024.100448}.

\bibitem[Sacilotto(2025)]{sacilotto2024riduzione}
Anna Sacilotto.
\newblock {Riduzione delle Emissioni di Carbonio mediante Ottimizzazione
  Stocastica dello Stoccaggio Energetico a Batteria in Sistemi Fotovoltaici}.
\newblock Master's thesis, University of Pavia, 2025.

\bibitem[Siler-Evans et~al.(2012)Siler-Evans, Azevedo, and
  Morgan]{silerevans2012marginal}
Kyle Siler-Evans, In{\^e}s~Lima Azevedo, and M.~Granger Morgan.
\newblock Marginal emissions factors for the {U.S.} electricity system.
\newblock \emph{Environmental Science \& Technology}, 46\penalty0 (9):\penalty0
  4742--4748, 2012.

\bibitem[Skoplaki and Palyvos(2009)]{skoplaki2009temperature}
E.~Skoplaki and J.~A. Palyvos.
\newblock On the temperature dependence of photovoltaic module electrical
  performance: A review of efficiency/power correlations.
\newblock \emph{Solar Energy}, 83:\penalty0 614--624, 2009.
\newblock \doi{10.1016/j.solener.2008.10.008}.

\bibitem[Tipping and Bishop(1999)]{tipping1999probabilistic}
Michael~E. Tipping and Christopher~M. Bishop.
\newblock Probabilistic principal component analysis.
\newblock \emph{Journal of the Royal Statistical Society: Series B (Statistical
  Methodology)}, 61\penalty0 (3):\penalty0 611--622, 1999.

\bibitem[van~der Meer et~al.(2021)van~der Meer, Wang, and
  Munkhammar]{vandermeer2021alternative}
Dennis van~der Meer, Guang~Chao Wang, and Joakim Munkhammar.
\newblock An alternative optimal strategy for stochastic model predictive
  control of a residential battery energy management system with solar
  photovoltaic.
\newblock \emph{Applied Energy}, 283:\penalty0 116289, 2021.
\newblock \doi{10.1016/j.apenergy.2020.116289}.

\bibitem[Wang et~al.(2014)Wang, Lin, and Pedram]{wang2014adaptive}
Yanzhi Wang, Xue Lin, and Massoud Pedram.
\newblock Adaptive control for energy storage systems in households with
  photovoltaic modules.
\newblock \emph{IEEE Transactions on Smart Grid}, 5\penalty0 (2):\penalty0
  992–1001, 2014.
\newblock \doi{10.1109/tsg.2013.2292518}.

\bibitem[Z{\"a}ngl et~al.(2015)Z{\"a}ngl, Reinert, R{\'i}podas, and
  Baldauf]{zangl2015icon}
G{\"u}nther Z{\"a}ngl, Daniel Reinert, Pilar R{\'i}podas, and Michael Baldauf.
\newblock The {ICON} ({ICOsahedral} non-hydrostatic) modelling framework of
  {DWD} and {MPI-M}: Description of the non-hydrostatic dynamical core.
\newblock \emph{Quarterly Journal of the Royal Meteorological Society},
  141:\penalty0 563--579, 2015.
\newblock \doi{10.1002/qj.2378}.

\bibitem[Zippenfenig(2024)]{zippenfenig2023openmeteo}
Patrick Zippenfenig.
\newblock Open-meteo.com weather {API}.
\newblock \url{https://doi.org/10.5281/zenodo.7970649}, 2024.

\end{thebibliography}

\end{document}